\documentclass[11pt]{article}
\usepackage{hyperref,amsmath,amsfonts,graphicx,psfrag,fancybox,qsymbols,indentfirst,amsthm,amssymb,amstext}
\usepackage[latin1]{inputenc}
\def \sous#1#2{\mathrel{\mathop{\kern 0pt#1}\limits_{#2}}}
\def \sur#1#2{\mathrel{\mathop{\kern 0pt#1}\limits^{#2}}}

\def \sp{^{(p)}}
\def \be{\begin{eqnarray*}}
\def \ee{\end{eqnarray*}}
\def \ben{\begin{eqnarray}}
\def \een{\end{eqnarray}}

\def\AArm{\fam0 }
\def\AAk#1#2{\setbox\AAbo=\hbox{#2}\AAdi=\wd\AAbo\kern#1\AAdi{}}%
\def\AAr#1#2#3{\setbox\AAbo=\hbox{#2}\AAdi=\ht\AAbo\raise#1\AAdi\hbox{#3}
}%

\def \sur#1#2{\mathrel{\mathop{\kern 0pt#1}\limits^{#2}}}
\def\BBp{{\AArm I\!P}}%

\def\E{{\math E}}

\def \e{{\rm e}}

\def \d{{\rm d}}

\def\E{\ensuremath{\mathbb{E}}}

\newcommand{\typ}{\operatorname{typ}}
\newcommand{\ac}{\operatorname{ac}}
\newenvironment{disarray}{\everymath{\displaystyle\everymath{}}\array} {\endarray}
\newcommand\nbOne{{\mathchoice {\rm 1\mskip-4mu l} {\rm 1\mskip-4mu l}
{\rm 1\mskip-4.5mu l} {\rm 1\mskip-5mu l}}}

\def\videbox{\mathbin{\vbox{\hrule\hbox{\vrule height1ex \kern.5em\vrule height1ex}\hrule}}}

\def \build#1#2#3{\mathrel{\mathop{\kern 0pt#1}\limits_{#2}^{#3}}}
\def\suPP#1#2{{\displaystyle\sup _{\scriptstyle #1\atop \scriptstyle #2}}}
\def\proDD#1#2{{\displaystyle\prod _{\scriptstyle #1\atop \scriptstyle #2}}}
\def\proof{\noindent{\bf Proof:}\hskip10pt}
\def\QED{\hfill\vrule height 1.5ex width 1.4ex depth -.1ex \vskip20pt}

\begin{document}
\newtheorem{lem}{Lemma}[section]
\newtheorem{defi}[lem]{Definition}
\newtheorem{theo}[lem]{Theorem}
\newtheorem{cor}[lem]{Corollary}
\newtheorem{prop}[lem]{Proposition}
\newtheorem{nota}[lem]{Notation}
\newtheorem{rem}[lem]{Remark}
\newdimen\AAdi%
\newbox\AAbo%
\font\AAFf=cmex10

\setcounter{page}{1}
\setcounter{section}{0}\date{}

\def \build#1#2#3{\mathrel{\mathop{\kern 0pt#1}\limits_{#2}^{#3}}}
\def\suPP#1#2{{\displaystyle\sup _{\scriptstyle #1\atop \scriptstyle #2}}}
\def\proDD#1#2{{\displaystyle\prod _{\scriptstyle #1\atop \scriptstyle #2}}}
\def\proof{\noindent{\bf Proof:}\hskip10pt}
\def\QED{\hfill\vrule height 1.5ex width 1.4ex depth -.1ex \vskip20pt}
\setcounter{page}{1}
\setcounter{section}{0}
\title{Martingales and  rates of presence in homogeneous fragmentations}
\large
\author{N. Krell\thanks {Université Rennes 1,
Campus de Beaulieu,
35042 Rennes Cedex. e-mail: nathalie.krell@univ-rennes1.fr}
\and
A. Rouault \thanks
{Université Versailles-Saint-Quentin, LMV, B\^atiment Fermat,
Universit\'e de Versailles 78035 Versailles Cedex. e-mail:
rouault@math.uvsq.fr}
}

\maketitle

\begin{abstract}
The main focus of 
this work is 
 the asymptotic behavior of mass-conservative homogeneous fragmentations. Considering the logarithm of masses makes the situation reminiscent of branching random walks. The standard approach is to study
  {\bf asymptotical} exponential rates (Berestycki \cite{Beres}, Bertoin and Rouault \cite{BertRou2}). For fixed $v > 0$, either the number of fragments whose sizes at time $t$ are of order $\e^{-vt}$  is exponentially growing  with rate $C(v) > 0$, i.e. the rate is effective, or the probability of presence of such fragments is exponentially decreasing with rate $C(v) < 0$, for some concave function $C$. In a recent paper \cite{krell}, N. Krell considered fragments whose sizes decrease at  {\bf exact} exponential rates, i.e. whose sizes are confined to be of order $\e^{-vs}$ for every $s \leq t$. In that setting, she characterized the effective rates. 
In the present paper we continue this analysis
 and focus on  probabilities of presence,   using the spine method and  a suitable martingale. 
For the sake of completeness, we compare our results   with those obtained in the standard 
 approach (\cite{Beres}, \cite{BertRou2}).

\end{abstract}


\noindent{\bf A.M.S. Classification. }  60J85, 65J25, 60G09.
\normalsize\rm
\bigskip

\section{Introduction and main results}
We begin with a brief overview on fragmentations. 
We refer the reader to Bertoin \cite{Bertoin1} for a more complete exposition (and also \cite{ba} and \cite{Beres}).
 We consider a homogenous fragmentation $F$ of intervals, which is a Markov process in continuous time taking its values in the set $\mathcal U$ of open sets of $(0,1)$. 
Informally, each interval component - or \textit{fragment} - splits as time goes on, independently of the others and with the same law, up to a rescaling.
We make  the restriction that the fragmentation is conservative, which means that no mass is lost.
  In this
case,   the law of $
 F $ is completely characterized by the so-called dislocation
measure $\nu$  (corresponding to the jump-component of the process)
which is a measure on $\mathcal{U}$ fulfilling the following conditions
$$\nu  ((0,1))=0,$$

  \begin{equation}\int_{\mathcal{U}}  (1-u_{1})\nu  (\d U)<\infty,
 \label{mesuredelevy}
 \end{equation}
and
$$\sum_{i=1}^{\infty}
u_{i}=1\ \ \ \ \ \hbox{for}\ \nu\!-\!\hbox{almost every} \ U\ \in \mathcal{U},$$  where for $U\in \mathcal{U}$,
$$|U|^{\downarrow}:=  (u_{1},u_{2},...)$$ is the decreasing
sequence of the  lengths of the interval components 
 of $U$.

It is a natural question  to study the rates of decay of fragments. If we measure the fragments by the logarithms of their sizes, a homogeneous fragmentation can be considered as an extension of 
a classical branching random walk in continuous time (\cite{Bertoin1} p. 21-22).  
For a broad range of branching models, the process either grows exponentially or becomes extinct.
 Let us recall some basic facts. 
If $\zeta_n$ is  a Galton-Watson process
 started from $\zeta_0 = 1$,  with finite mean $m = \mathbb E \zeta_1$, we have $n^{-1}\log \mathbb E(\zeta_n) = \log m$ and

(a)  \label{typea} if $m > 1$ and $\mathbb P(\zeta_1 \geq 1) = 1$, 
then a.s.\[\lim_n n^{-1}\log \zeta_n = \log m\ \ \ \ a.s.\]

(b) \label{typeb} if $m < 1$,  then   
\[\exists n_0 : \forall n \geq n_0 \ \ \ \zeta_n = 0\ \ \ \ a.s.,\]
 and 
\[\lim_n n^{-1} \log \mathbb P(\zeta_n \not= 0)= \log m\,.\]

More generally, in a branching random walk on $\mathbb R$, there is a concave function $\lambda$ which governs the local growth of the population.   
If $v$ is some speed and $Z_n$ is the number of particles 
 located around $nv$ in the $n$-th generation, then $\mathbb E Z_n$ grows exponentially at rate $\lambda(v)$. When $\lambda(v) > 0$, this quantity is also the effective exponential rate of growth of $Z_n$, (result of type (a), see \cite{biggins1977chernoff}). When $\lambda(v) < 0$, a.s. $Z_n$ is zero for $n$ large enough and $\lambda(v)$ is the effective exponential rate of decrease of $\mathbb P (Z_n \not=0)$ (result of type (b), see \cite{rouault1993precise}).

The goal of this paper is to present results of  the later kind (type (b)) for fragmentations, i.e. to study 
the asymptotic probability of 
 presence of abnormally large fragments.
Let us first explain known results of  type (a) - exponential growth -  and fix some notation. 

For $x\in (0,1)$ let $I_{x}  (t)$ be the  component of the
 interval fragmentation $F(t)$ which contains $x$, and let $|I_{x}(t)|$ be its
length.  Bertoin showed in \cite{ber2} that if $V$ is a uniform
random variable on $[0,1]$  independent of the
fragmentation, then $$ \xi (t):=-\log |I_{V}  (t)|$$
  is a subordinator  whose distribution is entirely determined by the characteristics of the
  fragmentation process. 
  Its Laplace transform is given by
\begin{equation}
\label{defLaplace}\mathbb E \e^{-q\xi(t)} = \e^{-t \kappa(q)} \ \ (q \geq 0)\end{equation}  
where $\kappa$ (the Laplace exponent) is the concave 
 function :
\ben
\label{defkappa}
\kappa  (q)
:=\int_{\mathcal{U}}\left  ( 1-\sum_{j=1}^{\infty}u_{j}^{q+1}\right)
\nu  (\d U)\,.
 \een
 In other words,
\[\kappa(q) = \int_{(0, \infty)} \left(1- \e^{-q x}\right) L(\d x)\,,\]
where  the Lévy measure $L$ is given by
\[L(dx) = \e^{-x} \sum_{j=1}^\infty \nu(-\log u_j \in \d x)\,.\]
If $\underline{p}$ is defined by
$$\underline{p}:=\inf\left \{p\in\mathbb{R}:\ \int_{\mathcal{U}}\sum_{j=2}^{\infty}u_{j}^{p+1}
\nu  (\d U)<\infty \right\}\,,$$
then Condition (\ref{mesuredelevy}) ensures that $\underline{p} \leq 0$ and we will assume throughout that $\underline{p} < 0$. It turns out that  $\kappa$ is an increasing concave analytic function on $(\underline{p}, \infty)$ and that (\ref{defLaplace}) holds also for $q\in (\underline{p} , 0)$. Set 
\[v_{\max} := \kappa'
 (\underline{p}^{+})\in [0, \infty]\,.\]
The SLLN tells us that  $\xi(t)/t \rightarrow \kappa'(0) =: v_{\typ}$ a.s., or in other words 
\[\lim_{t\rightarrow \infty}\  - t^{-1} \log |I_V(t)| = v_{\typ}\ \ \ \ a.s.\,.\]
To study effective exponential rates of decrease, we fix $a$ and $b$ two constants such that
\[0 < a <1 <b\,.\]
In the standard approach, one considers the set of fragments
\[\widetilde G_{v,a,b} (t) = \{I_x(t) : x \in (0,1) \ \ a e^{-vt} < |I_x(t)| < be^{-vt}\}\,.\]
To describe its behavior as $t \rightarrow \infty$, we need some notation. 
Define for $v < v_{\max}$
 \begin{equation}
\label{dual}
C(v) = \inf_{p> \underline p}\ (p+1) v - \kappa(p)\,.
\end{equation} 
This infimum is reached at a unique point $p= \Upsilon_v$ which is the unique solution to the equation $v = \kappa'(p)$, so that
\begin{equation}
\label{CV}
C(v) := (\Upsilon_v +1)v - \kappa(\Upsilon_v) 
\,.
\end{equation}
 If  $\overline p$ is the unique solution  of
the equation 
\[\kappa(q) =(q+1)\kappa'(q)\ \ \ q> \underline p\,.\]
 and if \[v_{\min} := \kappa'(\overline p)\,,\]
 then $C(v_{\min}) = 0$,   
 $C$ is
  positive for  $v \in (v_{\min}, v_{\max})$ and negative for $v < v_{\min}$. Moreover $C$ is concave analytic and $C'(v) = \Upsilon_v + 1$  (Legendre duality).

It is known  (\cite{Beres}, \cite{BertRou2}) that the asymptotic growth of $\widetilde G_{v,a,b} (t)$ 
is governed
 by $C(v)$ (which depends only on $v$ and not on $a,b$).
More precisely, we have\footnote{$\sharp A$ stands for the cardinality of the set $A$}:

$\bullet$ for $v\in (v_{\min}, v_{\max})$, then 
 a.s.
\begin{equation}
\label{BR0}
\lim_{t\rightarrow \infty}\ t^{-1} \log \sharp \widetilde G_{v,a,b}(t) = C(v)\; 
\end{equation}
 
$\bullet$  for $v < v_{\min}$, 
then  a.s.
\[\exists t_0 : \forall t \geq t_0 \ \ \ \sharp \widetilde G_{v,a,b}(t) = 0\,.
\]

 In a recent paper \cite{krell}, N. Krell  studied  the more constrained set  
\[G_{v,a,b}(t) = \{I_x(t) : x \in (0,1) \ \hbox{and}\ a e^{-vs} < |I_x(s)| < be^{-vs} \ \ \forall \!\ s \leq t\}\,,\]
and proved a result of the same kind. 
In particular, Proposition 3 (p.908) \cite{krell}   
tells us 
 that 
there exists 
 a positive number $\rho(v; a,b)$ depending 
 on  $v,a,b$ such that 

$\bullet$ for $v > \rho(v; a,b)$, 
 conditionally on $\{\inf\{t : G_{v,a,b}(t)=\emptyset\}=\infty\}$, a.s.
\begin{equation}
\label{krell0}
\lim_{t\rightarrow \infty}\ t^{-1} \log \sharp G_{v,a,b}(t) = v -\rho(v;a,b)\,, 
\end{equation}

$\bullet$ for $v < \rho(v; a,b)$, a.s.
\[\exists t_0 : \forall t \geq t_0 \ \ \ \sharp G_{v,a,b}(t) = 0\,.\]
Since the precise definition of $\rho(v;a,b)$ is rather involved, we postpone it in the forthcoming Section 3, formula (24). 

 This result  holds under the following assumption  A,  which 
 ensures the absolute continuity of the marginals of the underlying L\'evy process.
  Let $\nu_1$ be the pushforward of the measure $\nu$ by the mapping $U \mapsto u_1$.
\medskip

\noindent {\bf Assumption A}   
The absolutely continuous component $\nu_{1}^{\ac}$ of $\nu_1$ with respect to the Lebesgue measure on $[0,1]$ satisfies
\begin{equation}
\label{lac}
 \nu_{1}^{\ac}  ((1-\epsilon,1])=\infty \ \ \ \text{for any }\ \ \
\epsilon>0\,.
 \end{equation}

\medskip

The study of $\widetilde G_{v,a,b}$ or $G_{v,a,b}$ will be referred as the \textit{classical} or \textit{confined} model, respectively. 
According to the above informal classification of results 
 on branching models, we can say that
 the above assertions (\ref{BR0}) and (\ref{krell0}) are of type (a) on page \pageref{typea}. Our objective
  here is to 
present results of type (b) on page \pageref{typeb}.
 
For the \textit{classical} model,  an assumption is needed. A fragmentation is called $r$-lattice with $r >0$, if $\xi (t)$ is  a compound Poisson process whose jump measure has a support  carried by a discrete subgroup of $\mathbb R$ and $r$ is the mesh.
If there is no such $r$, the fragmentation is called non-lattice.
\medskip

\noindent {\bf Assumption B.} Either the fragmentation is non-lattice, 
or it is $r$-lattice and $a,b$ satisfy $b > a\e^r$.

\begin{theo}
\label{mtheo} \cite{BertRou3}
Under Assumption B, 
if $v <v_{\min} $, then
\begin{equation}
\label{vmintheo}
\lim_{t\rightarrow \infty}\ t^{-1} \log \mathbb P (\widetilde G_{v,a,b}(t) \not= \emptyset) = C(v)\,.
\end{equation}
\end{theo}
 In \cite{BertRou3}, the result of Theorem 5 is more precise since it gives sharp (i.e. non logarithmic) estimates of the latter probability.

For the more constrained set $G_{v,a,b}(t)$, 
 we have  the following theorem, which is the  
 main result of the present paper.
\begin{theo}
\label{wtheo}
Under  Assumption $A$, if $v -\rho(v; a,b) < 0$, then
\begin{equation}
\label{stheo}
\lim_{t\rightarrow \infty}\ t^{-1} \log \mathbb P (G_{v,a,b}(t) \not= \emptyset) = v-\rho(v; a,b)\,.
\end{equation}
\end{theo}
 Let us remark that since $G_{v,a,b}\subset \tilde{G}_{v,a,b}$, 
   the limits (\ref{krell0}) and (\ref{stheo}) are smaller than the limits (\ref{BR0}) and  (\ref{vmintheo}), respectively. 
   In fact we have the following general result, which extends  Remark 4 in  \cite{krell}. 
 \begin{prop}\label{cvrho} For all $v< v_{\max}$, then
  \ben C(v)\geq v-\rho(v; a,b).
  \een
  \end{prop}

 Let us explain shortly our method, whose  crucial tools were already 
 central in the proof of
  results of type (a) in \cite{krell}, namely 
 the construction 
 of  an additive martingale, 
 the corresponding change of probability and the 
 so-called spine decomposition.

 In branching or fragmentation problems, it is by now customary to enlarge the probability space by considering a randomly chosen branch or a randomly tagged fragment, respectively. This random element is called the "spine"\label{splinemethod}. Informally, we deal with two filtrations : a large one including the spine and a small one without the spine. A martingale built on the observation of the spine process may be projected on the small filtration, obtaining a so-called "additive" martingale. 
These martingales induce changes of probability. Making use of a proper choice of the martingale measurable with
 respect to the large
 filtration, the spine has generally a nice behavior under this new probability whereas the behavior of the other particles (or fragments) is not affected by
 this change.  
It is then possible to split the additive martingale into two parts : the contribution of the spine and the contribution of other terms (this is called the spine decomposition). It allows to describe in the large time limit 
 the behavior of the additive martingale itself and the behavior of the branching or fragmentation. 

For the \textit{classical} model, the Esscher martingale is convenient to study $\widetilde G_{v,a,b}$ (see \cite{BertRou2}).  For the branching random walk and related processes, a good recent reference with historical comments is \cite{hardy2009spine}.
The change of probability forces the tagged fragment to have a "good" asymptotic logarithmic rate of decrease.

  For the \textit{constrained} model, 
 as in \cite{krell}, we have been inspired by  the change of probability introduced by Bertoin \cite{ber1} and Lambert \cite{lam}.  It has the effect of forcing the spine to be confined in some interval, as required to study  $G_{v,a,b}$.

In Section 2, we summarize the basic notions on fragmentations and Lévy processes which will be needed later.  
  Section 3 is devoted to the study of the two martingales and their asymptotic properties.   
  In Section 4, we give the proofs of the theorems on the presence probabilities and the proofs of the results on martingales\footnote{In particular,  a mistake in the proof of Theorem 2.1 in \cite{krell} is corrected.}. 
   For the sake of completeness,  a direct short proof of 
Theorem \ref{mtheo} with the spine method  is given to illustrate the common feature of both models
  (it should be noted that a similar method 
was applied to obtain analogous results in the context of   
    branching Brownian motion in
  \cite{HardyHarrisPres}).
Section 5 is devoted to a proof  of Proposition \ref{cvrho}, only based on properties of L\'evy processes.
 
\section{Background on fragmentations and  L\'evy processes.}

\subsection{Partition fragmentations and interval fragmentations} \label{subsect21}
Let $\mathbb N$ stand for the set of positive integers; a block is a subset of $\mathbb N$. For every $k \in \mathbb N$, the block
$\{1,...,k\}$ is denoted by $[k]$.
Let $ \mathcal{P}$ the space of
partitions of $ \mathbb{N}$.
 As in \cite{berrou}, 
 any
measure:
$$\omega =\sum_{ (t,\pi,k)\in\mathcal{D}}^{\infty} \delta_{ (t,\pi,k)},$$
where $\mathcal{D}$ is a subset of
$\mathbb{R}_{+}\times\mathcal{P}\times\mathbb{N}$ such that
\begin{equation}
\label{beres}
\forall t^{'}\geq 0\ \ \forall n\in\mathbb{N}\ \
\sharp\left\{ (t,\pi,k)\in\mathcal{D} \ |\ t\leq t^{'},
\pi_{|[n]}\neq ([n],\emptyset,\emptyset,...),k\leq n\right\}<\infty
\end{equation}
and for all $t\in\mathbb{R}$
\[\omega(\{t\}\times \mathcal{P}\times \mathbb{N})\in\{0,1\}.\]
is called a  discrete point measure on the space
$\Omega:=\mathbb{R}_{+}\times\mathcal{P}\times\mathbb{N}$.
Starting from an arbitrary discrete point measure $\omega$ on
$\mathbb{R}_{+}\times\mathcal{P}\times\mathbb{N}$, we will construct
a nested partition $\Pi= (\Pi (t),t\geq 0)$ (which means that for
all $t\geq t^{'}$ $\Pi (t)$ is a finer partition of $ \mathbb{N}$
than $\Pi (t^{'})$). We fix $n\in\mathbb{N}$, the assumption (\ref{beres}) that
the point measure $\omega$ is discrete enables us to construct a
step path $ (\Pi (t,n),t\geq 0)$ with values in the space of
partitions of $[n]$, which only jumps at times $t$ at which the
fiber $\{t\}\times\mathcal{P}\times\mathbb{N}$ carries an atom of
$\omega$, say $ (t,\pi,k)$, such that $\pi_{|[n]}\neq
([n],\emptyset,\emptyset,...)$ and $k\leq n$. In that case, $\Pi
(t,n)$ is the partition obtained by replacing the $k$-th block of
$\Pi (t-,n)$, denoted $\Pi_{k} (t-,n)$, by the restriction
$\pi_{|\Pi_{k} (t-,n)}$ of $\pi$ to this block, and leaving the
other blocks unchanged. Of course for all $t\geq 0$, $ (\Pi
(t,n),n\geq 0)$ is compatible (i.e. for every $n$, $\Pi (n,t)$ is a
partition of $[n]$ such that the restriction of $\Pi (n+1,t)$ to
$[n]$ coincide with $\Pi (n,t)$). As a consequence, there exists a
unique partition $\Pi (t)$, such that for all $n\geq 0$ we have $\Pi
(t)_{|[n]}=\Pi (t,n)$. 
This process $\Pi$ is a
 partition-valued homogeneous fragmentation (\cite{Bertoin1} chap. 3).

Let the set $\mathcal{S^{\downarrow}}$ be $$\mathcal{S^{\downarrow}}:=\left\{s= (s_{1},s_{2},...)\ |\
 s_{1}\geq s_{2}\geq ...\geq 0\ , \sum_{i=1}^{\infty}
s_{i}\leq 1 \right\}.$$ 

\noindent A block $B$ 
 has an asymptotic
frequency, if the limit
\[|B|:=\lim_{n\rightarrow\infty}n^{-1} \sharp(B\cap [n])\] exists.
 When every block of some partition $\pi\in \mathcal{P}$
 has an asymptotic frequency, we write
 $|\pi|= (|\pi_{1}|,...)$ and then $|\pi|^{\downarrow}= (|\pi_{1}|^{\downarrow},...)\in\mathcal{S}^{\downarrow}$ for the
 decreasing rearrangement of the sequence $|\pi|$.
 If a block of the partition $\pi$ does
 not have an asymptotic frequency, we decide that
 $|\pi|=|\pi|^{\downarrow}=\partial$, where $\partial$
 stands for some extra point added to $ \mathcal{S}^{\downarrow}$.

On $\Omega$,  the sigma-field generated 
  by the
restriction  to $[0,t]\times\mathcal{P}\times\mathbb{N}$ is denoted by $
\mathcal{G} (t)$. So $ {\mathcal G} = (\mathcal{G} (t))_{t\geq 0}$ is a
filtration,  and the nested partitions $ (\Pi (t),t\geq 0)$ are 
${\mathcal G}$-adapted. If $|\Pi (r)|^{\downarrow}$ is the
decreasing rearrangement  of the sequence of
the asymptotic frequencies of the blocks of $\Pi (r)$, we denote by $ \mathcal{F} (t)$ 
  the
sigma-field   generated by $\left(|\Pi (r)|^{\downarrow}\right)_{r \leq t}$. 
Of course ${\mathcal F} =(\mathcal{F} (t))_{t\geq 0}$ is a sub-filtration of $\mathcal G$.

 Let
$\mathcal{G}_{1} (t)$ the sigma-field generated by the
restriction of the discrete point measure $\omega$ to the fiber
$[0,t]\times\mathcal{P}\times\{1\}$. So ${\mathcal G}_1 = \left({\mathcal G}_1(t)\right)_{t\geq 0}$ 
 is a sub-filtration of ${\mathcal G}$, 
 and the
first block of $\Pi$ is ${\mathcal G}_1$-measurable.
Let $\mathcal{D}_{1}\subseteq \mathbb{R}_{+}$ be the random set of
times $r\geq 0$ for which the discrete point measure has an atom on
the fiber $\{r\}\times\mathcal{P}\times\{1\}$, and for every
$r\in\mathcal{D}_{1}$, denote the second component of this atom by
$\pi (r)$.

 There is a powerful link between partition fragmentations and interval fragmentations. On the one hand, the $\mathcal{S^{\downarrow}}$-valued process of ranked asymptotic frequencies  $|\Pi|^{\downarrow}$ of a partition fragmentation is a so-called ranked (or mass) fragmentation (\cite{be2}, \cite{ber2}), and conversely a partition fragmentation 
can be built from a ranked fragmentation via a "paintbox" process. On the other hand, 
the interval decomposition $(J_{i}  (t), J_{2}  (t),...)$
of the open $F  (t)$ ranked in decreasing order is a ranked fragmentation, denoted by $ X(t):= (|J_{i}  (t)| , |J_{2}  (t)| , ... )^{\downarrow}$. 
We can then lift this ranked fragmentation to a partition fragmentation. More precisely, if $\nu$ is the dislocation measure of $F$, and $\tilde \nu$ its image by the map $U \mapsto |U|^{\downarrow}$, 
then
  according to
Theorem 2 in \cite{ber2}, there exists a unique measure $\mu$ on $
\mathcal{P}$ which is exchangeable (i.e. invariant by the action of finite permutations on $ \mathcal{P}$), and such that $\tilde\nu$ is the
image of $\mu$ by the map  $\pi\mapsto |\pi|^{\downarrow}$ 
where
 $|\pi|^{\downarrow}$ is  the decreasing
rearrangement of the sequence of the asymptotic
frequencies of the blocks of $\pi$. 
So,  for all measurable function $f:
[0,1]\rightarrow\mathbb{R}_{+}$  such that $f (0)=0$,
\[\int_{\mathcal{P}}f (|\pi_{1}|)\mu (\d\pi)= \int_{\mathcal{S^{\downarrow}}}\sum_{i=1}^{\infty}s_{i}f (s_{i})\tilde\nu (\d s)= \int_{\mathcal{U}}\sum_{i=1}^{\infty}u_{i}f (u_{i})\nu (\d U)\,.\]
It should be noted that  $\{|J_{1}(t)|,|J_{2}(t)|,... \}_{t\geq 0} = \{|\Pi_{1}(t)|,|\Pi_{2}(t)|,...\}_{t\geq 0}$ (equality in distribution in general, and true equality if $\Pi$ is obtained by a paintbox process). 

In the following sections, $\Pi$   refers to this partition fragmentation.

\subsection{L\'evy processes.}\label{subsectionlevy}
A L\'evy process is a stochastic process with c\`adl\`ag sample paths and stationary independent increments (\cite{Bertoinb96}). In this work, two types of such processes will play a key role :  

$\bullet$ a subordinator is a L\'evy process taking values in $[0,\infty)$, which implies that its sample paths are increasing,

$\bullet$ a L\'evy process is
  completely asymmetric  when 
its jumps are either all positive or all negative.
   We will consider here L\'evy processes without positive jumps, i.e. spectrally negative processes. 

The Laplace transform of a subordinator $\sigma = (\sigma_t)_{t\geq 0}$ is given by 
\footnote{Bold symbols $\mathbf{P}$ and
$\mathbf{E}$ will refer to L\'evy processes while $\mathbb{P}$ and
$\mathbb{E}$ refer to  fragmentations.
} :
\ben
\label{defLaplacesub}
{\mathbf E} \exp -\lambda \sigma_t = \exp -t \Phi(\lambda) , \ \ \ \lambda \geq 0\,,\een
where $\Phi$ is called the Laplace exponent. If ${\mathcal E} = ({\mathcal E}_t)_{t\geq 0}$ 
is the natural filtration associated with $\sigma$ 
\[\left(\exp (-p \sigma_t + t \Phi(p))\right)_{t \geq 0}
\]
is a ${\mathcal E}$-martingale. We define the probability measure $\mathbf{P}^{(p)}$ as the 
Esscher transform:
\begin{equation}
\label{change0}
\d\mathbf{P}^{(p)} |_{\mathcal{E}_{t}} = \exp\{-p \sigma_t + t \Phi(p)\}\ \d\mathbf{P} |_{\mathcal{E}_{t}}\,.
\end{equation} 
Under $\mathbf{P}^{(p)}$, 
 $\sigma$ is a subordinator with Laplace exponent $q \mapsto \Phi(p+q) - \Phi(p)$. 
 The change of probability forces the process to have mean $t\Phi'(p)$ at time $t$. It also exponentially tilts the Lévy measure. 
 
Let us recall some facts about 
 about
completely asymmetric L\'evy processes, 
lifted  from \cite{Bertoinb96} and \cite{ber1}. Let $Y=
( Y_{t})_{t\geq0}$ be a L\'evy process with no positive jumps and let ${\mathcal E}= 
  (\mathcal{E}_{t})_{t\geq 0}$  be the natural filtration associated with $Y$. 
 The case where $Y$ is the negative of a
subordinator
 is degenerate for our purpose and
is  therefore implicitly excluded 
 in the rest of the paper.
The law of the L\'evy process started at $x\in \mathbb{R}$ will be
denoted by $\mathbf{P}_{x}$, its Laplace transform is
given by
$$\mathbf{E}_{0}  (\e^{\lambda Y_{t}})= \e^{t\psi  (\lambda)}, \ \ \lambda \geq0,$$ where $\psi:
\mathbb{R}_{+}\rightarrow\mathbb{R}$ is called the Laplace exponent. The function $\psi$ is convex with
$\lim_{\lambda\rightarrow\infty}\psi  (\lambda)=\infty$.

Let $\phi:\mathbb{R}_{+}\rightarrow\mathbb{R}_{+}$ be the right
inverse of $\psi$  
so that $\psi  (
\phi   (\lambda))=\lambda$ for every $\lambda\geq0$.
The scale function $W:\mathbb{R}_{+}\rightarrow\mathbb{R}_{+}$ 
 is  the unique
 continuous function  with Laplace transform:
$$\int_{0}^{\infty}e^{-\lambda x} W  (x) \d x=\frac{1}{\psi  (\lambda)}\
\ \ , \ \lambda>\phi  (0).$$
For $q\in\mathbb{R}$, let $W^{   (
q)}:\mathbb{R}_{+}\rightarrow\mathbb{R}$ be the
 continuous function defined by 
\[W^{  (q)}  (x):=\sum_{k=0}^{\infty}q^{k}W^{*k+1}  (x)\,,\]
 where $W^{*k}$ is  the $k$-th convolution of
 $W$ with himself,  
 so that
\[\int_{0}^{\infty}e^{-\lambda x} W^{ (q)}  (x) \d x=\frac{1}{\psi  (\lambda)-q}\
\ \ , \ \lambda>\phi  (q)\,.\]
For fixed $x$, $W^{(q)}(x)$ can be seen as an analytical function in $q$.
The functions $W^{(q)}$ are useful to investigate the two-sided exit problem for  L\'evy processes.  Their properties are well exposed in the book of Kyprianou \cite{Kypbook} and in \cite{KypChan}, examples are in \cite{KypHub} and in \cite{kyprianou13special}.

 The following theorems \ref{theo11} and \ref{ergodic} taken from \cite{lam}
and  \cite{ber1}  yield another important  martingale and its corresponding change of probability.

\begin{theo}\label{theo11}  Suppose that the one-dimensional distributions of the L\'evy
 process $Y$ are absolutely continuous.
 Let us define  the
 critical value
 \begin{equation}
  \rho_{\beta}:=\inf\{ q\geq 0\ ;\
 W^{  (-q)}  (\beta)=0\}\label{rho}.
 \end{equation}
 Then the following holds:
 \begin{enumerate}
\item $\rho_{\beta}\in  (0,\infty )$
 and the function $W^{  (-\rho_{\beta})}$ is strictly positive on $ (0,\beta )$.
 \item For $\beta >0$, let 
\begin{equation}
T_{\beta}=\inf\{t : \ Y_{t} \notin  (0,\beta )\}\ ; \label{ta}
 \end{equation} then the 
process $D = (D_t)_{t\geq 0}$ with 
\begin{equation}
D_{t}:=e^{\rho_{\beta} t}\ \mathbf{1}_{\{t<T_{\beta}\}}\ \frac{W^{
  (-\rho_{\beta})}  (Y_{t})}{ W^{  (-\rho_{\beta})}  (x)}\label{dt}
\end{equation}
 is a $  (\mathbf{P}_{x},  \mathcal{E})$-martingale,  for every $x \in (0, \beta)$.
 \item The mapping
$\beta\mapsto\rho_{\beta}$ 
 is
strictly decreasing and continuous  on $(0,\infty)$. 
 \end{enumerate}
\end{theo}

Point 1 is taken from \cite{ber1} Theorem 2(i) and (iii). Point 2 is from \cite{lam} Theorem 3.1 (ii). Point 3 is from \cite{lam} Proposition 5.1 (ii). 
\begin{rem}
\label{rem2}
\begin{enumerate}
\item
Notice that Proposition 5.1 in \cite{lam}, devoted to the smoothness properties of functions $W$ and $\rho$,
  is claimed for paths with unbounded variation or with bounded variation provided the Lévy measure has no atoms. 
However, 
 this additional assumption is not used in the part of the proof dedicated to our point 3, it is used to prove stronger regularity.
 However we do not need this assumption, since we only care about continuity. 
  Let us stress that the smoothness of scale functions is a very active subject, see \cite{KypChan}. 
\item
 The definition of $\rho_{\beta}$ is  complicated, but some examples are given in \cite{krell}.
\end{enumerate}
 \end{rem}

 Let the probability measure $\mathbf{P}^{\updownarrow}$ be the $h$-transform of $\mathbf{P}$ based on
 the martingale $D$ :
\begin{equation}\d\mathbf{P}_{x}^{\updownarrow} |_{\mathcal{E}_{t}}\
=D_{t}\!\ \d\mathbf{P}_{x}
|_{\mathcal{E}_{t}}.\label{changementproba3}\end{equation}

 \begin{theo}
\label{ergodic}
With the same assumption as in Theorem \ref{theo11}, 
under $\mathbf{P}^{\updownarrow}_x$, 
$Y$ is a homogeneous strong Markov process on $(0,\beta )$, positive-recurrent and as $t\rightarrow\infty$, $Y_{t}$  converges in distribution to its stationary probability, which  has a density.
\end{theo}
This is essentially a rephrasing of  Theorem 3.1 in \cite{lam}, the convergence in distribution is a consequence  of Theorem 2 (v) of \cite{ber1}.

The change of probability forces the process to be confined in $(0,\beta)$. In the Brownian motion case, this is sometimes called 
 a taboo process \cite{Knight}.

\section{Two additive martingales and their asymptotic behavior}
One of the  most striking fact about homogeneous fragmentations is the subordinator representation. If $\xi_t = -\ln |\Pi_1(t)|$, then,  
as seen in \cite{Bertoin1} (Section 3.2.2), the process $\xi = \left(\xi_t\right)_{t\geq 0}$ defined on $(\Omega, {\mathcal G}, \mathbb P)$ is a  subordinator, which means in particular that $\xi_{t+s} -\xi_t$ is independent of  ${\mathcal G}(t)$. In this section, we will adapt the  statements of Section \ref{subsectionlevy} to the subordinator $\xi$ (instead of $\sigma$) and to the spectrally negative process $Y = \left(Y_t = vt - \xi_t -\log a\right)_{t \geq 0}$, starting at $x= -\log a$. It should be stressed that  ${\mathcal G}$ is not the proper filtration of these processes, but the martingale properties remain true, as well as the Markov property. We will then perform a projection on the filtration ${\mathcal F}$ of the ranked fragmentation. 

\subsection{The classical additive martingale $M\sp$}
\label{classical}
As in (\ref{change0}), we define for $p>\underline{p}$ the $\mathcal{G}$-martingale 
$D^{(p)} = \left(D_{t}^{(p)}\right)_{t\geq 0}$ where \[D_{t}^{(p)}= \e^{-p\xi (t) +t \kappa (p)}\,,\]  
and the probability measure   $\mathbb{P}^{(p)}$ as the  transform :
\begin{equation}\d\mathbb{P}^{(p)} |_{\mathcal G (t)}\
=D_{t}^{(p)} \d\mathbb{P}
|_{\mathcal G (t)}.\label{changementproba2}\end{equation}

Projecting the martingale $D^{(p)}$ 
 on the sub-filtration ${\mathcal F}$,
 we obtain  the well-known additive 
${\mathcal F}$-martingale $M^{(p)} = \left(M_{t}^{(p)}\right)_{t\geq 0}$, where
\begin{equation}
  \label{add2}
  M_{t}^{(p)}=\sum_{j=1}^{\infty}|\Pi_{j}(t)|^{p+1}e^{\kappa (p)t}=\sum_{i=1}^{\infty}|J_{i}(t)|^{p+1}e^{\kappa (p)t}\,.\end{equation}
 The projection of  
 (\ref{changementproba2}) 

gives the identity:
\begin{equation}
\label{change3}
\d\mathbb{P}^{(p)} |_{\mathcal{F}(t)}\
=M_{t}^{(p)} \d\mathbb{P} |_{\mathcal{F}(t)}\,.\end{equation}
In \cite{berrou} Proposition 6, there is a complete description of the behavior of the process $\Pi$. We keep in mind the next result.

\begin{lem}\label{lemmexip} Under $\mathbb{P}^{(p)}$, the process $\xi$ is a subordinator with Laplace exponent $q \mapsto \kappa (p+q)-\kappa (p)$.
\end{lem}

\subsection{The martingale $M^{(v,a,b)}$ and the confined process.} 
\label{associated}

Since we are interested in the set of the \textquotedblleft good intervals\textquotedblright~ at time $t$ as
\begin{equation}
G_{v,a,b}(t)=\left\{I_{x}  (t): \ x\in (0,1) \ \ \hbox{and} \  \ a e^{-vs}<
|I_{x}  ( s)|< b e^{-vs} \ \ \forall\ s\leq t\right\}\label{good}
\end{equation}
it is natural to study the Lévy process $Y = (Y_t)_{t\geq 0}$ defined by 
\[Y_t := vt -\xi_t-\log a\]
and its exit time from $(0 , \log b/a)$.
Clearly $Y$
 has no positive
jump and its
 Laplace exponent is 
\begin{equation}
\label{3.2}\psi  (\lambda) = v\lambda-\kappa  (
\lambda)\,,\end{equation} with $\kappa$ defined in (\ref{defkappa}).

Assumption (A) guarantees  that the marginals of 
 $Y$ are absolutely continuous and  Theorem~\ref{theo11} can thus be applied.
Let us fix some notation. The distribution of $Y$ depends on $v$ and we set
\begin{equation}
\label{Tetrho}
\rho(v; a,b):=\rho_{\log  (b/ a)} \ , \ T:=T_{\log  (b/a)}\,
\end{equation}
where $\rho_{\beta}$ and  $T_{\beta}$ are defined in (\ref{rho}) and (\ref{ta}), respectively. 
We will use frequently the notation $\rho$ instead of $\rho(v; a,b)$. 

 To further simplify the notation, let also
\begin{equation}\label{defh}h  (y):= W^{  (-\rho)}  (y-\log
a)\mathbf{1}_{\{y\in  (\log a ,\log b )\}}
\end{equation} for all
$y\in\mathbb{R},$ and $h (-\infty)=0$.
This function is well defined thanks to
Theorem~\ref{theo11} and  $h(0) \not= 0$.

By 
translating  (\ref{dt}) 
into the new notation  we get a $
\mathcal{G}$-martingale $D= (D_t)_{t\geq 0}$ 
\ben\label{defDh}D_{t}= \e^{\rho t}\ \mathbf{1}_{\{t<T\}}\ \frac{h
  (vt-\xi_t)}{h  (0)}\,,\een
  and then a new probability defined by
  \begin{equation}
\d\mathbb{P}^{\updownarrow} |_{\mathcal{G}(t)}\
=D_{t}\!\ \d\mathbb{P}
|_{\mathcal{G}(t)}\,.\label{changementproba}
\end{equation}

For $i\geq 1$, let $P_{i} (t)$ be the block of $\Pi (t)$ which contains $i$ at time
$t$.  We define the killed partition as follows \[\Pi^{\dag}_{j}
(t)=\Pi_{j} (t)\mathbf{1}_{\{\exists i \in \mathbb{N}^{\ast}| \
\Pi_{j} (t)=P_{i} (t);\ \forall s\leq t\ \  |P_{i} (s)|\in
(ae^{-vs},be^{-vs})\}}.\]

Similarly, if $I$ is an interval component of $F  (t)$, we define the
``killed'' interval $I^{\dag}$ by $I^{\dag}=I$ if $I$ is good  (i.e.
$I\in G_{v,a,b}
 (t)$ with $ G_{v,a,b}  ( t)$ defined in~ (\ref{good})), else by
$I^{\dag}=\emptyset$. Projecting 
 the martingale $D$ on the
sub-filtration $\mathcal{F}$, 
 we obtain  an
additive martingale $M^{(v,a,b)} = \left(M_{t}^{(v,a,b)}\right)_{t\geq 0}$ where
\begin{eqnarray*}\nonumber M_{t}^{(v,a,b)}&=&\frac{e^{\rho t}}{h  (0)} \sum_{j\in\mathbb{N}}
h\left  (vt+\log |\Pi_{j}^{\dag}  ( t)|\right)|\Pi_{j}^{\dag}
 (t)|\\ &=&\frac{e^{\rho t}}{h  (0)} \sum_{i\in\mathbb{N}}
h\left  (vt+\log |J_{i}^{\dag}  ( t)|\right)|J_{i}^{\dag}
 (t)|\,.
\end{eqnarray*}

Finally, let  the  absorption time of $M^{(v,a,b)}$ at $0$ be
 $$\begin{disarray}{rcl}
\zeta:=\inf\{t :M_{t}^{(v,a,b)} =0\}=\inf\{t: G_{v,a,b}
(t)=\emptyset\},\end{disarray}$$ with the convention
$\inf{\emptyset}=\infty$.

The projection of (\ref{changementproba}) on ${\mathcal F}$ 
 gives the identity:
\begin{equation}
\label{change4}\d\mathbb{P}^{\updownarrow} |_{\mathcal{F}(t)}\
=M_{t}^{(v,a,b)} \d\mathbb{P} |_{\mathcal{F}(t)}\,.\end{equation}

 The
 upshot is that the change of probability $\mathbb P^{\updownarrow}$ only affects the behavior of the block 
 which contains 1. More precisely, like in  lemma 8  (ii) \cite{BertRou2}, we obtain:
\begin{lem}\label{lm1}
Suppose Assumption (A) holds. Under $\mathbb{P}^{\updownarrow}$, the
 restriction of $\omega$ to
 $\mathbb{R}_{+}\times\mathcal{P}\times\{2,3,...\}$ has the same
 distribution as under $\mathbb{P}$\,. 
\end{lem}
\subsection{Growth of martingales}
The next 
 theorems govern the asymptotic behaviors of our martingales $M^{(v,a,b)}$ and $M\sp$, according to the values of parameters $v$ and $p$.

 It should be noted that assertion 1 of Theorem \ref{theo2} was claimed in \cite{krell} Theorem 2, but unfortunately there was a mistake in the proof. 
Indeed it is not true in general that the function $h$ defined in (\ref{defh}) is Lipschitz at $0$. 
See again Remark \ref{rem2} for smoothness of $W$ (and hence of $h$). 

The points 1) and 2 a) of Theorem \ref{cvmtp} are known (\cite{Bertoin03} p. 406-407 and \cite{berrou} respectively). We will recall the argument for the sake of completeness.
\begin{theo}\label{theo2}
Suppose assumption A holds, then:
 \begin{enumerate}
\item
If $v>\rho(v; a,b)$, the martingale $M^{(v,a,b)}$ is bounded in $\mathrm{L}^{2}  (\mathbb{P})$.
\item If $v < \rho(v; a,b)$,   

a) $\lim_{t\rightarrow\infty} M^{(v,a,b)}_{t} = 0$, $\mathbb P$-a.s.,

b) there exists $K_1, K_2 > 0$ such that for every $t$
\begin{equation}
\label{mtsouscr}
K_1 \leq \e^{(v -\rho(v; a,b))t}\!\ \mathbb E \left[M^{(v,a,b)}_t\right]^2 \leq K_2\,.
\end{equation}
 \end{enumerate}
\end{theo}

\begin{theo}
\label{cvmtp}
\begin{enumerate}\item If 
$p\in (\underline p, \overline p)$,  there exists $\alpha > 0$ such that the martingale $M\sp$ is bounded in $\mathrm{L}^{1+\alpha}  (\mathbb{P})$.
\item If  $p \geq \overline p$,

a) $\lim_{t\rightarrow \infty} M_t\sp = 0$, $\mathbb P$-a.s.

b) There exists $\alpha_0 > 0$ such that for $\alpha \in (0, \alpha_0)$, 
\[d(p, \alpha) := (1+\alpha)\kappa(p) -\kappa\left((1+\alpha) (p+1) -1\right) > 0\]
and then for those $\alpha$, we have for  every $t > 0$ 
\begin{equation}
\label{esppositive}
\e^{d(p,\alpha) t}\!\ \mathbb E\!\ |\! M_t^{(p)}\! |^{1+\alpha} \leq C_{\alpha,p}\,,
\end{equation}
where  $C_{\alpha, p}$ depends  on $\alpha$ and $p$. 
\end{enumerate}
\end{theo}  

\section{Proofs}
\subsection{Proof of Theorem \ref{wtheo}}
\proof $\bullet$ We first show the upper bound of (\ref{stheo}), i.e.
\begin{equation}
\label{limsup}
\limsup_{t\rightarrow \infty}\ t^{-1} \log \BBp (G_{v,a,b}(t) \not= \emptyset) \leq v-\rho(v; a,b)\,.
\end{equation}

 Let $0< \bar{a}<a<1<b<\bar{b} $. As in Section \ref{associated}, we associate with $\bar a, \bar b$ and $v$, the parameter $\bar{\rho} =\rho(v; \bar a, \bar b)$, as well as the set of "good" intervals 
$$\bar G(t)= G_{v, \bar a, \bar b}  (t):=\{I_{x}  (t) : \ \  x\in (0,1) \ \  \hbox{and}\ \  \ |I_{x}  (s)|\ \in \
 (\bar{a} e^{-vs},\bar{b} e^{-vs} )\ \ \ \forall\ s\leq t\}\,,$$
and the martingale   
 $\bar{M} = M^{v, \bar{a}, \bar{b}}$.

Let for every $y\in\mathbb{R}$
\[\bar{h}  (y):=
W^{  (-\bar{\rho})}  (y -\log  \bar{a})\mathbf{1}_{\{y\in  (\log
\bar{a},\log \bar{b})\}}\,.\]
 For $t\geq 0$ fixed, we have:

\be
1=\mathbb{E}\bar{M}_{t}&=& \frac{ e^{\bar{\rho} t}}{\bar{h}   (
0)}\mathbb{E} \left(\sum_{i\in\mathbb{N}} \bar{h}  (vt+\log |J_{i}  (t)|) |J_{i}
  (t)|\ \mathbf{1}_{\{ J_{i}  (t)\in \bar G  (t)\}}\right)
\\&\geq &
 \frac{\bar{a} e^{(\bar{\rho}-v)
t}}{\bar{h}  (0)}\mathbb{E}\left(\sum_{i\in\mathbb{N}} \bar{h}  (vt+\log |J_{i}
  (t)|)\
   \ \mathbf{1}_{\{J_{i}  (t)\in G_{v,a,b}  (t)\}}\right)\,.
\ee
Since $ (a, b)\subsetneq  (\bar{a},
\bar{b})$,  the function $\bar h$ is continuous and strictly positive on $[ \log a, \log b ]$ 
so that, if
$$K_{3}:=\bar{h}  (0)/\left   (
\bar{a}\underset{x\in[\log a, \log b]}{\inf}\bar{h}   (
x)\right)<\infty,$$ 
then,  for all $t\geq 0$ :
\[K_{3}\geq e^{  (\bar{\rho}-v)t}\mathbb{E}\left(\sum_{i\in\mathbb{N}}
\mathbf{1}_{\{ J_{i}  (t)\in G_{v,a,b}  (t)\}}\right)\geq  e^{  (\bar{\rho}-v)t}\mathbb{P} 
(G_{v,a,b}(t)\not=\emptyset)\,,\]
and consequently
\[\limsup_{t\rightarrow\infty}\ t^{-1} \log\mathbb{P}
(G_{v,a,b}(t)\not=0)\leq v-\bar{\rho} = v - \rho(v; \bar a, \bar b)\,.\]
Since it holds true for all $\bar{a}, \bar{b}$ such that $0< \bar{a}<a<1<b<\bar{b}$, we can let
 $\bar{a}\rightarrow a$ and $\bar{b}\rightarrow b$  and use the
continuity of $\rho(v; \cdot , \cdot)$ (see Theorem \ref{theo11}.3) to 
obtain  the inequality  (\ref{limsup}).

$\bullet$
It remains to prove the lower bound of  (\ref{stheo}), i.e.
\begin{equation}
\label{liminf}
\liminf_{t\rightarrow\infty}\ t^{-1} \log \BBp(G_{v,a,b}(t) \not= \emptyset) \geq v-\rho\,.
\end{equation}
Let us 
drop the subscripts and superscripts  $(v,a,b)$.
Since $M$ is a positive martingale and 
  $\{G(t) \not = \emptyset\} = \{M_{t} \not = 0\}  $, we have
 $$1=\mathbb{E}(M_{t})=\mathbb{E}(M_{t}\nbOne_{\{M_{t}\not= 0\}}) = \mathbb{E}(M_{t}\nbOne_{\{G(t)\not= \emptyset\}})\,.$$ 
Now, thanks to  the Cauchy-Schwarz inequality:
\[
\mathbb{E}(M_{t}\nbOne_{\{G(t)\not= \emptyset\}})\leq \left(\mathbb{E}(M_{t}^{2})\right)^{1/2}\left(\mathbb{P}(G(t)\not=\emptyset)\right)^{1/2}\,\] and applying  (\ref{mtsouscr}), we get 
\[\mathbb{P}(G(t)\not= \emptyset)\geq K_2^{-1}  e^{
(v-\rho) t}\,,\]
which yields (\ref{liminf}).
\QED

\subsection{Proof of Theorem \ref{mtheo}}
The upper bound is straightforward. For all $p\geq \underline p$, we have \[1=\mathbb{E} M_{t}^{(p)}=\mathbb{E}\left(\sum_{i=1}^{\infty}|J_{i}(t)|^{p+1}\e^{\kappa (p)t}\right)\geq a^{p+1}\e^{\kappa(p) t -(p+1) vt} \mathbb{P}(\widetilde G_{v,a,b} (t) \not= \emptyset)\,.\]

Hence \[\mathbb{P}(\widetilde G_{v,a,b} (t) \not= \emptyset)\leq a^{-(p+1)}\e^{[(p+1)v-\kappa (p)]t}\,\]
and
\[\limsup_{t\rightarrow\infty}\ t^{-1}\log \mathbb{P}(\widetilde G_{v,a,b} (t) \not= \emptyset)\leq (p+1)v-\kappa (p).\]
 In particular,  for $p=\Upsilon_{v}$, we get from (\ref{CV})
\[\limsup_{t\rightarrow\infty}\ t^{-1}\log \mathbb{P}(\widetilde G_{v,a,b} (t) \not= \emptyset)\leq C(v)\,.\]
To prove the lower bound
\begin{equation}
\label{liminftilde}
\liminf_{t\rightarrow\infty}\ t^{-1}\log \mathbb{P}(\widetilde G_{v,a,b} (t) \not= \emptyset)\geq C(v)\,,\end{equation}
we   use  again the  change of probability  (\ref{change3}) with $p= \Upsilon_v$.
We have, 
\[\mathbb P(\widetilde G_{v,a,b} (t )\not= \emptyset) = \mathbb E\sp\left((M_t\sp)^{-1} ; \widetilde G_{v,a,b} (t )\not= \emptyset\right)\geq \ \ \ \ \ \ \ \ \ \ \ \ \ \ \ \ \ \ \ \ \ \ \ \ \ \ \ \ \ \ \ \ \ \ \ \ \ \ \ \ \ \ \ \ \ \]
\begin{equation}\label{mino1} \ \ \ \ \ \ \ \ \ \ \ \ \ \ \ \ \ \e^{t C(v)-t\varepsilon}\!\  \mathbb P\sp\left(\sup_{0 < s \leq t} M_s\sp < \e^{-t C(v)+t\varepsilon}; vt-\xi_t  
\in[\log a, \log b]\right)  
\end{equation}
and
\[\mathbb P\sp\left(\sup_{0 < s \leq t} M_s\sp < \e^{-t C(v)+t\varepsilon} ; vt -\xi_t \in[\log a, \log b]\right) \geq\ \ \ \ \ \ \ \ \ \ \ \ \ \ \ \ \ \ \ \ \ \ \ \ \ \ \ \ \ \ \ \ \ \ \ \ \ \ \ \ \ \ \ \ \ \]
\begin{equation}
\label{mino2}
{\mathbb P}\sp\left( vt - \xi_t \in[\log a, \log b]\right) - \mathbb P\sp\left(\sup_{0 < s \leq t} M_s\sp \geq \e^{-t C(v)+t\varepsilon}\right)\,.
\end{equation}

From Lemma \ref{lemmexip} we see that under ${\mathbb P}\sp$, the L\'evy process  
$\left(vt -\xi_t\right)_{t \geq 0}$ has mean $-\kappa'(p) +v=0$ and variance 
$\sigma^2_p := -\kappa''(p)$. From Proposition 2 of Bertoin and Doney \cite{BertIHP}
it satisfies the local central limit theorem, if it is not lattice. We get
\begin{equation}
\label{lclt}
\sigma_p\sqrt{2\pi t}\ \mathbb P\sp (vt - \xi_t \in [\log a, \log b]) \rightarrow \log\frac{b}{a}\,.
\end{equation}
and then
\begin{equation}
\label{lcltcor}
\liminf_t \ t^{-1} \log \mathbb P\sp (vt - \xi_t \in [\log a, \log b]) =0\,.
\end{equation}
In the case of a $r$-lattice fragmentation, under assumption B, there is at least an integer multiple of $r$ in the interval
 $[vt-\log b, vt -\log a]$. We can use the lattice version of the local central limit theorem (see for instance \cite{Chat} Theorem 2 iii)),  to obtain (\ref{lcltcor}) again.

To tackle the second term of the RHS of (\ref{mino2}), we argue as in  \cite{HardyHarrisPres}.
 By convexity $((M_t\sp)^{1+\alpha}, t \geq 0)$ is a $\mathbb P$-submartingale, so $((M_t\sp)_{t\geq 0}^\alpha , t \geq 0)$ is a $\mathbb P\sp$-submartingale. 
Hence, by Doob's inequality,  
\begin{eqnarray}\nonumber
\mathbb P\sp\left(\sup_{0 < s \leq t}\mid\!M_s\sp\!\mid\!\ \geq \e^{-t C(v)+t\varepsilon}\right)&\leq&  \e^{t\alpha C(v)-\alpha t\varepsilon}\E\sp \mid\! M_t\sp\!\mid^\alpha\\ \label{Doob}
&=& \e^{t\alpha C(v)-\alpha t\varepsilon} \E \mid\!M_t\sp\!\mid^{1+\alpha}\,,
\end{eqnarray}
and by (\ref{esppositive})
 we have for $\alpha \in (0, \alpha_1]$ for some $\alpha_1 > 0$
\begin{eqnarray}
\label{Doob2}\mathbb P\sp\left(\sup_{0 < s \leq t}\mid\!M_s\sp\!\mid\!\ \geq \e^{-t C(v)+t\varepsilon}\right)
\leq K'_{\alpha, p} \e^{tH(\alpha)}\,,
\end{eqnarray}
 where
\[ H(\alpha) = \alpha C(v) - \alpha\varepsilon + 
d(p, \alpha)
\,,\]
and $K'_{\alpha, p}$ is some constant.
Now, a second order development of $\kappa$ around $p$ gives
\begin{eqnarray*} H(\alpha) 
= - \alpha\varepsilon - \frac{\alpha^2(p+1)^{2}}{2}\kappa''(p)(1+o(\alpha))
\end{eqnarray*}
and, since $\kappa'' < 0$ ($\kappa$ is concave),  we may choose $\alpha$  small enough such that $H(\alpha) < 0$. This yields 
\begin{equation}
\label{gr}
\limsup_t \ t^{-1}\log \mathbb P^{(p)} \left(\sup_{0 < s \leq t}\mid\!M_s\sp\!\mid\!\ \geq \e^{-t C(v)+t\varepsilon} \right) < 0.
\end{equation}

So, gathering (\ref{gr}) and (\ref{lclt}), 
 we obtain
\[\liminf_{t\rightarrow \infty} t^{-1} \log \mathbb P\sp\left(\sup_{0 < s \leq t} M_s\sp < \e^{-t C(v)+t\varepsilon}; vt - \xi_t \in[\log a, \log b]\right) = 0\,,\]
which, with (\ref{mino1}), yields 
\[\liminf_{t\rightarrow \infty} t^{-1} \log \mathbb P(\widetilde G_{v,a,b} (t) \not= \emptyset)\geq C(v) - \varepsilon\]
for all $\varepsilon > 0$. Letting $\varepsilon\rightarrow 0$ proves (\ref{liminftilde}), and concludes the proof of Theorem \ref{mtheo}.
\QED
\subsection{Proof of Theorem \ref{theo2} :}
We use  the change of probability (\ref{change4}): 
\begin{equation}
\label{m2}\mathbb{E
}(M_{t}^{2})=\mathbb{E}^{\updownarrow}(M_{t})\,,\end{equation} and the spine decomposition (see page \pageref{splinemethod} for a discussion of this method): 
$$M_{t}=c_{t}+d_{t},$$
where
\begin{equation}\label{defct}c_{t}:=\frac{e^{\rho t}}{h (0)}h\left (vt+\log (|\Pi_{1}^{\dag} (t)|)\right)
\ |\Pi_{1}^{\dag} (t)|\end{equation} and
\begin{equation}\label{defdt}d_{t}:=\frac{e^{\rho t}}{h (0)} \sum_{i=2 }^{\infty}
h\left (vt+\log\left (|\Pi^{\dag}_{i} (t)|\right)\right)\
|\Pi^{\dag}_{i} (t)|\,.
\end{equation}
The asymptotic behaviors of $c_t$ and $d_t$  are governed by the two following lemmas.
\begin{lem}
\label{a2c} 
Suppose Assumption (A) holds. 
Under $\mathbb P^{\updownarrow}$,  $\e^{-(\rho-v)t}c_t$ converges in distribution as $t \rightarrow \infty$, to a random variable $\eta$ with no mass at $0$. Moreover  there exists $K >0$ such that 
\begin{equation}
\label{e2c} 
\lim_{t\rightarrow \infty} \mathbb{E}^{\updownarrow} (c_{t}) \e^{-(\rho-v) t} = K\,. 
\end{equation}
\end{lem}
\begin{lem} 
\label{a2d}Suppose Assumption (A) holds. If $\rho\not= v$, there exists $L > 0$ such that
\begin{equation}
\label{e2d}
\mathbb E^{\updownarrow} d_t \leq L \max\{\e^{(\rho -v)t}, 1\}\,.
\end{equation}
\end{lem}

\subsubsection{Proof of Theorem \ref{theo2} 1)} 
 From (\ref{m2}),
 it is enough to
prove that
$\lim_{t\rightarrow \infty}\mathbb{E}^{\updownarrow} (M_{t})<\infty$. Now, 
by (\ref{e2c}), we have $\lim_{t\rightarrow \infty}
\mathbb{E}^{\updownarrow} (c_{t})=0$ and  
by (\ref{e2d}), we have $\sup_t \mathbb E^{\updownarrow} (d_{t}) < \infty$. \QED

\subsubsection{Proof of Theorem \ref{theo2} 2) a)}
\label{432}  The method is now classic, (see for instance  \cite{Kypalt}) 
and uses a decomposition which may be found e.g. in  Durrett \cite{Dur} p. 241. It will be used also  in the proof of Theorem \ref{cvmtp} below.  We  only need to prove that  $\mathbb P^\updownarrow(\limsup M_t = \infty) =1$.

We have the obvious lower bound
$$M_t \geq c_t$$
For $v <\rho$,  Lemma (\ref{a2c})  yields
 $\lim c_t = \infty$ in $\mathbb P^\updownarrow$ probability, or in other words 
$\mathbb P^\updownarrow(\limsup_t c_t = \infty) =1$ 
 which  implies $\mathbb P^\updownarrow(\limsup_t M_t = \infty) =1$, hence $\mathbb P(\lim_t M_t = 0) =1$.\QED

\subsubsection{Proof of Theorem \ref{theo2} 2 b)} 
It is a straightforward consequence of (\ref{m2}) and Lemma \ref{a2c} and \ref{a2d}. \QED

\subsubsection{Proof of Lemma \ref{a2c} :}
 From the definition (\ref{defct}) of $c_t$,  we see that the distribution under $\mathbb P^{\updownarrow}$ of the process $\left(\e^{-(\rho-v)t}c_t  ,  t\geq 0\right)$  is the same as the distribution under 
$\mathbf{P}_{\log (1/a)}^{\updownarrow}$ of
$\left(h (0)^{-1} W^{(-\rho)} (Y_{t})\e^{Y_{t}}\ \mathbf{1}_{\{t<T\}}  ,  t\geq 0\right)$. Under the  latter probability,
 the stopping time $T$ (defined in (\ref{Tetrho})) is a.s. infinite and  from Theorem 
\ref{ergodic}, $Y$ is positive-recurrent, it
 converges in distribution and the limit has no mass in $0$.   
Since the function  $y \mapsto W^{(-\rho)}(y)\e^y$ is continuous, it is bounded on the compact support of the distribution of $Y_t$,
and  $y_t = \mathbb{E}^{\updownarrow} \left(c_t \e^{-(\rho -v )t}\right)$ has a positive limit.
\QED
 
\subsubsection{Proof of Lemma \ref{a2d} :}
We start from the definition of $d_t$ decomposing the time interval $[0,t]$ into pieces $[k-1, k[$ and splitting the sum (\ref{defdt}) according to the time where the fragment separates from $1$:
\begin{equation}
\label{hed}
h(0)\e^{-\rho t} d_t = \sum_k S_k
\end{equation}
with
\[S_k(t) = \sum_{i \in {\mathcal I}_k} h(vt + \log |\Pi^{\dag}_{i} (t)|)|\Pi^{\dag}_{i} (t)|\]
where ${\mathcal I}_k$ is the set of $i \geq 2$ such that the block $\Pi_i(t)$  separates at some instant $r \in {\mathcal D}_1 \cap[k-1, k[$. The block after the split which contains $1$ is $\Pi_1(r)$. Thus, there is some index $\ell \geq 2$ such that $\Pi_i(t)  \subseteq \pi_\ell(r) \cap\Pi_1(r-)$. Then, at time $k$, $\pi_\ell(r) \cap\Pi_1(r-)$ is partitioned into $\Pi_j(k), j \in {\mathcal J}_{\ell, r}$ where ${\mathcal J}_{\ell,r}$ is some set of indices measurable with respect to
${\mathcal G}(k).$ Conditionally  upon 
${\mathcal G}(k)$, the partition $(\Pi_i(t) , i \in {\mathcal I}_k)$ can be written in the form
$\tilde\Pi^{(j)} (t-k)_{|\Pi_{j} (k)} , j \in {\mathcal J}_k$, where  

$\bullet$ ${\mathcal J}_k$ is some set of indices ${\mathcal G}(k)$-measurable 

$\bullet$ $(\tilde\Pi^{(j)})_{j\in\mathbb{N}}$
is a family of i.i.d. homogeneous fragmentations distributed as $\Pi$ under
$\mathbb{P}$ (see Lemma \ref{lm1}).

As a consequence:
\begin{equation}\label{underset}\underset{i\in {\mathcal I}_k}{\cup}\Pi_{i} (t)=\cup_{j\in {\mathcal J}_k}\tilde\Pi^{(j)} (t-k)_{|\Pi_{j} (k)}\,,\end{equation}
with the slight abuse of notation by which we write a union of partitions instead of the union of the blocks of these partitions,
 and  for all $m\in\mathbb N$ \begin{equation}\label{tilde}|\tilde{\Pi}_{m}^{(j)} (t-k)_{|\Pi_{j} (k)}|
 =|\tilde{\Pi}_{m}^{(j)} (t-k)||\Pi_{j} (k)|.\end{equation}
Now, we have to take into account the killings.

Let us call \textquotedblleft good fragment\textquotedblright a fragment which   satisfies the constraint of non-killing all along its history up to time $t$. We can decompose
\be S_k (t)=
 \sum_{j \in {\mathcal J}_k}|\Pi_{j} (k)|\left(\sum_m h\left(vt + \log(|\tilde{\Pi}_{m}^{(j)} (t-k)||\Pi_{j} (k)|)\right)|\tilde{\Pi}_{m}^{(j)} (t-k)| \nbOne_{j,m,k}\right)\ee
where $\nbOne_{j,m,k} = 1$ if and only if  $\tilde{\Pi}_{m}^{(j)} (t-k)_{|\Pi_{j} (k)}$ is a good fragment. If $\Pi_{j} (k)$ is a good fragment, we define
\[\widetilde M^{j}_{t-k} := \e^{-\rho(k-t)}\sum_m \frac{h(vt + \log(|\tilde{\Pi}_{m}^{(j)} (t-k)||\Pi_{j} (k)|))}{h(vk + \log(|\Pi_{j} (k)|)}|\tilde{\Pi}_{m}^{(j)} (t-k)| \nbOne_{j,m,k}\,.\]
From the  definition of $h$ in (\ref{defh}), the process $\left(\widetilde M^{j}_{t-k}\right)_{t\geq k}$ is a 
$({\mathcal G}(s))_{s\geq k}$ martingale with respect to $\mathbb P$, distributed as $(M_t)_{t\geq 0}.$ 

Denoting $\nbOne_{j,k} = 1$ if and only if   $\Pi_{j} (k)$ is a good fragment, from Lemma \ref{lm1} we have  
\be\mathbb{E}^{\updownarrow} \left( S_k(t)| {\mathcal G}(k)\right) = \e^{\rho(k-t)} \sum_{j \in {\mathcal J}_k} |\Pi_{j} (k)| h(vk + \log(|\Pi_{j} (k)|) \nbOne_{j,k}\,.\ee
Now again by the definition of $h$ and its  continuity, there exists $C_{3}>0$ such that 
\[ h(vk + \log(|\Pi_{j} (k)|) \leq C_{3} \nbOne_{vk + \log(|\Pi_{j} (k)| \in (\log a, \log b)}\]
and
\[\mathbb{E}^{\updownarrow} \left( S_k(t)| {\mathcal G}(k)\right) \leq C_{3} \e^{(\rho-v)k}\e^{-\rho t}\sum_{j \in {\mathcal J}_k} \nbOne_{j,k}\,.\]
It is clear that
the only terms that contribute to the above sum 
 correspond to  good fragments at time $k$ 
 which were dislocated from good $\Pi_1(k-1)$ during $[k-1, k[$.  
Since the fragmentation is conservative, there were  at most  $be^{v}/a$ such dislocations 
 during that time, which yields:
 \[\mathbb{E}^{\updownarrow} \left( S_k(t)| {\mathcal G}(k)\right) \leq C_{3} be^{v}a^{-1}\e^{(\rho-v)k}\e^{-\rho t}\,.\]
Coming back to (\ref{hed}) we get 
 \be\mathbb{E}^{\updownarrow} (d_{t})\leq 
 C_{4}\sum_{k=1}^{\lfloor t\rfloor}\e^{(\rho-v)k}\,,\ee
for some constant $C_{4}>0$. 
In other words, for all $v\not=\rho$ there exists $L>0$ such that
\[\mathbb{E}^{\updownarrow} (d_{t})\leq L \max (\e^{
(\rho-v) t},1)\,,\]
which proves (\ref{e2d}), hence Lemma \ref{a2d}.\QED 
 
\subsection{Proof of Theorem \ref{cvmtp}}
Let us recall the definition of 
 the function 
\[d(p, \alpha) = (1+\alpha) \kappa(p) - \kappa \big((1+\alpha)(p+1) - 1\big)\]
and let us look at its sign. We have $d(p,0)=0$ and the derivative of $d(p,\alpha)$ in $\alpha=0$ is 
\begin{equation}
\label{cases}
\kappa(p) - (p+1)\kappa'(p)\
\begin{cases} &\!< 0\;\;\;\;\mbox{if}\; p < \bar p,\\
 &\!=0\;\;\;\;\mbox{if}\; p = \bar p,\\
      &\!>0\;\;\;\;\mbox{if}\; p > \bar p\,.
\end{cases}
\end{equation} 
If $p = \bar p$, the second derivative in $\alpha=0$ is $-(\bar p +1) \kappa''(\bar p) >0$ since $\kappa$ is concave.
This ensures that, in any case, there exists $\alpha_0(p) >0$ such that, for every $\alpha \in (0, \alpha_0(p))$
\begin{equation}
\label{casesbis}
d(p,\alpha)\
\begin{cases} &\!< 0\;\;\;\;\mbox{if}\; p < \bar p,\\
      &\!>0\;\;\;\;\mbox{if}\; p \geq \bar p\,.
\end{cases}
\end{equation} 

 1) 
 In the proof of Theorem 2 of \cite{Bertoin03} p.406-407, Bertoin gave the

 estimate:
\begin{equation}
\label{bdg}
\E \sup_{0 < s \leq t}\mid\! M_s\sp\!\mid^{1+\alpha}\!\ \leq 
K_\alpha\!\ c(p,\alpha)\!\ \int_{0} ^{t} \exp\left(d(p, \alpha)s
\right) \!\ \d s
\end{equation}
where $K_\alpha$ is a universal constant depending only on $\alpha$, 
and
\[c(p,\alpha) = \int_{{\mathcal S}^*}\mid\sum_{i=1}^\infty (x_i^{p+1} - x_i)\mid^{1+\alpha}\!\ \nu(\d x) < \infty\]
for  every $p > \underline p$ and $\alpha\in [0, \alpha_1(p)]$ for some $\alpha_1(p) >0$. 
From (\ref{casesbis}) above,  the integral on the RHS of (\ref{bdg}) is then uniformly bounded in $t$.
 
2) a) 
The martingale is  bounded by below by the contribution of the spine: 
$$M_t\sp \geq \e^{t\kappa(p)}|\Pi_1(t)|^{p+1} = \exp \{t\kappa(p) - (p+1) \xi_t\}\,.$$
As an easy consequence of Lemma \ref{lemmexip} , we see that under $\mathbb{P}\sp$, 
 the L\'evy process $\big(\kappa(p)t - (p+1) \xi_t\big)_{t\geq 0}$ has mean $\kappa(p) - (p+1)\kappa'(p)$ 
 which is nonegative from (\ref{cases}) since $p \geq \bar p$. 
We get successively
$\mathbb P\sp(\limsup_{t\rightarrow \infty} (\kappa(p)t - (p+1) \xi_t) = \infty)=1$, hence
  $\mathbb P\sp(\limsup M_t\sp = \infty)=1$, and  $\mathbb P(\lim M_t\sp = 0)$ (see section \ref{432}). 

 2 b) The only point which remains to prove is  (\ref{esppositive}), and it is a consequence of (\ref{bdg}) and (\ref{casesbis}).

\section{Comparison of  limits.}
\noindent{\it Proof of Proposition \ref{cvrho}:}
Let us give a direct proof of the inequality
\begin{equation}\label{v-rho}v- \rho(v; a,b) \leq C(v)\,.\end{equation} 

Fix $v,a,b$ and let $\rho= \rho(v; a,b)$ and $\beta= \log b/a$. 
By the definition (\ref{dual}) of $C$, (\ref{v-rho}) is equivalent 
to
\begin{equation}\label{interm}pv -\kappa(p) \geq -\rho\end{equation} 
for every $p > \underline{p}$. Referring to  
(\ref{3.2}) we recall  that $pv -\kappa(p)= \psi(p)$, 
(the Laplace exponent  of the process $Y$), 
so that (\ref{interm}) is equivalent to
\begin{equation}
\label{last}\rho+\psi(p) \geq 0\,.\end{equation}
If $\psi(p) \geq 0$, there is nothing to prove since $\rho$ is nonegative by definition. 
 Let us assume  $\psi(p) < 0$.
 If $W$ is the scale function of $Y$, we have 
\[
\rho =  \inf \{q \geq 0 : W^{(-q)}(\beta) =0\} 
=  \inf \{q' \geq \psi(p) : W^{(\psi(p)-q')}(\beta) =0\} -  \psi(p)\,.
\]
Moreover, if $W_p$ is the scale function of the tilted process of Laplace exponent $\lambda \mapsto \psi(\lambda +p) - \psi(p)$, 
we have
\[W^{(\psi(p)-q')}(\beta) = e^{px}W_p^{(-q')}(x)\]
(\cite{Kypbook} p.213 and Lemma 8.4 p.222).
 This yields
\[\rho+ \psi(p) =  \inf \{q' \geq \psi(p) : W_{p}^{(-q')}(\beta) =0\}\,.\]
Since $W_{p}^{(-q')}(\beta) > 0$ for $q' \leq 0$, the latter infimum is nonegative, so that  (\ref{last}) holds true, 
which ends the proof.
 \QED

\begin{rem} A consequence of this proposition is that  when $v<v_{min}$, we have $\rho(v; a,b)>v$.
\end{rem}
\newpage

\textit{Acknowledgements.} We would like to thank an Associate Editor and two anonymous referees
for their pertinent comments and suggestions on the first draft.
\bibliographystyle{plain}
\small
\bibliography{bib0705}

\begin{thebibliography}{10}

\bibitem{ba}
A-L. Basdevant.
\newblock Fragmentation of ordered partitions and intervals.
\newblock {\em Electron. J. Probab.}, 11:no. 16, 394--417, 2006.

\bibitem{be2}
J.~Berestycki.
\newblock Ranked fragmentations.
\newblock {\em ESAIM Probab. Statist.}, 6:157--175 (electronic), 2002.

\bibitem{Beres}
J.~Berestycki.
\newblock Multifractal spectra of fragmentation processes.
\newblock {\em J. Statist. Phys.}, 113(3-4):411--430, 2003.

\bibitem{Bertoinb96}
J.~Bertoin.
\newblock {\em L\'evy processes}, volume 121 of {\em Cambridge Tracts in
  Mathematics}.
\newblock Cambridge University Press, Cambridge, 1996.

\bibitem{ber1}
J.~Bertoin.
\newblock Exponential decay and ergodicity of completely asymmetric {L}\'evy
  processes in a finite interval.
\newblock {\em Ann. Appl. Probab.}, 7(1):156--169, 1997.

\bibitem{ber2}
J.~Bertoin.
\newblock Homogeneous fragmentation processes.
\newblock {\em Probab. Theory Related Fields}, 121(3):301--318, 2001.

\bibitem{Bertoin03}
J.~Bertoin.
\newblock The asymptotic behavior of fragmentation processes.
\newblock {\em J. Europ. Math. Soc.}, 5(4):395--416, 2003.

\bibitem{Bertoin1}
J.~Bertoin.
\newblock {\em Random fragmentation and coagulation processes}, volume 102 of
  {\em Cambridge Studies in Advanced Mathematics}.
\newblock Cambridge University Press, Cambridge, 2006.

\bibitem{BertIHP}
J.~Bertoin and R.~A. Doney.
\newblock Spitzer's condition for random walks and {L}\'evy processes.
\newblock {\em Ann. Inst. H. Poincar\'e Probab. Statist.}, 33(2):167--178,
  1997.

\bibitem{berrou}
J.~Bertoin and A.~Rouault.
\newblock Additive martingales and probability tilting for homogeneous
  fragmentations.
\newblock Report 808 of University Paris 6 available at\\
  \href{http://www.proba.jussieu.fr/mathdoc/textes/PMA-808.pdf}{http://www.pro%
ba.jussieu.fr/mathdoc/textes/PMA-808.pdf}, 2003.

\bibitem{BertRou3}
J.~Bertoin and A.~Rouault.
\newblock Asymptotical behaviour of the presence probability in branching
  random walks and fragmentations.
\newblock
  \href{http://fr.arxiv.org/PS_cache/math/pdf/0409/0409547.pdf}{ar{X}iv:math.P%
R/04095477}, 2004.

\bibitem{BertRou2}
J.~Bertoin and A.~Rouault.
\newblock Discretization methods for homogeneous fragmentations.
\newblock {\em J. London Math. Soc. (2)}, 72(1):91--109, 2005.

\bibitem{biggins1977chernoff}
J.D. Biggins.
\newblock {Chernoff's theorem in the branching random walk}.
\newblock {\em J. Appl. Probab.}, 14(3):630--636, 1977.

\bibitem{KypChan}
T.~Chan, A.E. Kyprianou, and V.~Savov.
\newblock Smoothness of scale functions for spectrally negative {L}\'evy
  processes.
\newblock
  \href{http://arxiv.org/PS_cache/arxiv/pdf/0903/0903.1467v1.pdf}{ar{X}iv/0903%
.1467}, 2009.

\bibitem{Chat}
S.D. Chatterji.
\newblock Asymptotic formulae derived from the central limit theorem.
\newblock {\em Confer. Sem. Mat. Univ. Bari}, 234:37, 1990.

\bibitem{Dur}
R.~Durrett.
\newblock {\em Probability: theory and examples}.
\newblock Duxbury Press, Belmont, CA, second edition, 1996.

\bibitem{hardy2009spine}
R.~Hardy and S.~Harris.
\newblock {A Spine Approach to Branching Diffusions with Applications to
  $L^p$-Convergence of Martingales}.
\newblock In {\em S{\'e}minaire de Probabilit{\'e}s XLII}, volume 1979 of {\em
  Lecture Notes in Math.}, pages 281--330. Springer, Berlin, 2009.

\bibitem{HardyHarrisPres}
R.~Hardy and S.C. Harris.
\newblock A conceptual approach to a path result for branching {B}rownian
  motion.
\newblock {\em Stochastic Process. Appl.}, 116(12):1992--2013, 2006.

\bibitem{KypHub}
F.~Hubalek and A.E. Kyprianou.
\newblock Old and new examples of scale functions for spectrally negative
  {L}\'evy processes.
\newblock
  \href{http://arxiv.org/PS_cache/arxiv/pdf/0801/0801.0393v2.pdf}{ar{X}iv/0801%
.0393}, 2008.

\bibitem{Knight}
F.B. Knight.
\newblock Brownian local times and taboo processes.
\newblock {\em Trans. Amer. Math. Soc.}, 143:173--185, 1969.

\bibitem{krell}
N.~Krell.
\newblock Multifractal spectra and precise rates of decay in homogeneous
  fragmentations.
\newblock {\em Stochastic Process. Appl.}, 118:897--916, 2008.

\bibitem{Kypalt}
A.E. Kyprianou.
\newblock Travelling wave solutions to the {K}-{P}-{P} equation: alternatives
  to {S}imon {H}arris' probabilistic analysis.
\newblock {\em Ann. Inst. H. Poincar\'e Probab. Statist.}, 40(1):53--72, 2004.

\bibitem{Kypbook}
A.E. Kyprianou.
\newblock {\em Introductory lectures on fluctuations of {L}\'evy processes with
  applications}.
\newblock Universitext. Springer-Verlag, Berlin, 2006.

\bibitem{kyprianou13special}
A.E. Kyprianou and V.~Rivero.
\newblock {Special, conjugate and complete scale functions for spectrally
  negative L{\'e}vy processes}.
\newblock {\em Electron. J. Probab.}, 13:1672--1701, 2008.

\bibitem{lam}
A.~Lambert.
\newblock Completely asymmetric {L}\'evy processes confined in a finite
  interval.
\newblock {\em Ann. Inst. H. Poincar\'e Probab. Statist.}, 36(2):251--274,
  2000.

\bibitem{rouault1993precise}
A.~Rouault.
\newblock {Precise estimates of presence probabilities in the branching random
  walk}.
\newblock {\em Stochastic Process. Appl.}, 44(1):27--39, 1993.

\end{thebibliography}
\end{document}